\documentclass{article}
\usepackage{amsmath}
\usepackage{amsfonts}
\usepackage[top=1.5in,bottom=1.5in,left=1.5in,right=1.5in]{geometry}

\numberwithin{equation}{section}
\newtheorem{thm}{Theorem}[section]

\newtheorem{lem}[thm]{Lemma}

\newtheorem{ex}{Example}[section]

\newcommand{\be}{\begin{equation}}
\newcommand{\ee}{\end{equation}}
\newcommand{\ben}{\begin{enumerate}}
\newcommand{\een}{\end{enumerate}}

\newcommand{\qed}{\hspace*{\fill}Q.E.D.}  
\title{\Large On the classification of projectively flat Finsler metrics with constant flag curvature}
\author{Benling Li\footnote{Research is support in part by  NNSFC(10801080), ZPNSFC(LY13A010013) and K.C. Wong Magna Fund in Ningbo University. } }

\date{August 11, 2013}

\begin{document}

\maketitle
\begin{abstract}
In this paper, we study locally projectively flat Finsler metrics with constant flag curvature ${\bf K}$.
We prove those are totally determined by their behaviors at the origin by solving some nonlinear PDEs.
The classifications when ${\bf K}=0$, ${\bf K}=-1$ and ${\bf K} =1$ are given respectively in
an algebraic way.
Further, we construct a new projectively flat Finsler metric with flag curvature ${\bf K}=1$ determined by a Minkowskian norm with double square roots at the origin. As an application of our main theorems, we give the classification of locally projectively flat  spherical symmetric Finsler metrics much easier than before.
\end{abstract}
{\bf Keywords}: Finsler metric; projectively flat; constant flag curvature

\section{Introduction}

    The regular case of Hilbert's Fourth Problem is to study and characterize Finsler metrics on an open subset in $R^n$ whose geodesics are straight lines.
Such metrics are called {\it locally projectively flat } Finsler metrics. Riemannian metrics form a special and important class in Finsler geometry.
Beltrami's theorem tells us that a Riemannian metric is locally projectively flat if and only if it is with constant sectional curvature ${\bf K} = \lambda$, which can be expressed as
\be \label{Rie}
F_{\lambda} = \frac{ \sqrt{|y|^2 + \lambda ( |x|^2|y|^2 - \langle x, y \rangle^2}) }{1+\lambda|x|^2},
\ee
where $ y\in T_x\mathcal{U} \approx R^n$, $\mathcal{U}\subset R^n$.
However, it is not true in general.

   Flag curvature is an analogue of sectional curvature in Finsler geometry.  It is known that there are many locally projectively flat Finsler metrics which are not with constant flag curvature; and there are many Finsler metrics with constant flag curvature which are not locally projectively flat. A natural problem is to characterize projectively flat Finsler metrics with constant flag curvature.
   In \cite{F1}\cite{F2}, P. Funk classified projectively flat Finsler metrics with constant flag curvature on convex domains in $R^2$.
The famous Funk metric $ F = F(x,y) $ defined on unit ball ${\rm B}^n$ in $R^n$ is locally projectively flat with flag curvature ${\bf K} = -\frac{1}{4}$ is given by
 \be \label{Funk}
 F = \frac{ \sqrt{(1-|x|^2)|y|^2 + \langle x, y \rangle^2}}{1-|x|^2} + \frac{\langle x, y \rangle}{1-|x|^2},
 \ee
 where $y\in T_x {\rm B}^n \approx R^n$.
In 1929, L. Berwald studied locally projectively flat Finsler metrics, specially in the case of zero flag curvature \cite{Be1}\cite{Be2}.
He gave the equivalent equations of such metrics and found that the key problem is to solve the following PDE:
\be \label{main}
\Phi_{x^k} = \Phi \Phi_{y^k},
\ee
where $\Phi = \Phi(x,y)$, $x$, $y \in R^n$. However, it is difficult to solve above equation at that time
though he constructed a projectively flat Finsler metric with ${\bf K}=0$ which be called Berwald's metric now  as following
\be \label{Berwaldmetric}
B= \frac{ (\sqrt{ (1-|x|^2) |y|^2 + \langle x, y \rangle^2} +\langle x, y \rangle )^2}{ (1-|x|^2)^2 \sqrt{ (1-|x|^2) |y|^2 + \langle x, y \rangle^2} },
\ee
where $ y\in T_x{\rm B}^n \approx R^n$.
The first locally projectively flat non-Riemannian Finsler metric with positive flag curvature $K=1$ was given by R. Bryant on $S^2$ \cite{3}. By algebraic equations, Z. Shen  gave the
following expression of Bryant's example including the higher
dimension in \cite{Sh1}.
\begin{equation}\label{Bryant}
\begin{split}
 F(x,y)&= \mathcal{I}m \Big[ \frac{-\langle x,y\rangle + i \sqrt{(e^{2i\alpha} + |x|^2)|y|^2 - \langle x,y\rangle^2 }}{e^{2i\alpha} + |x|^2}\Big]\\
 &=\sqrt{\frac{\sqrt{\mathcal {A}}+\mathcal
{B}}{2\mathcal {D}}+(\frac{\mathcal {C}}{\mathcal
{D}})^2}+\frac{\mathcal {C}}{\mathcal {D}},
\end{split}
\end{equation}
where
\begin{eqnarray}
\nonumber&&\mathcal {A}:=(|y|^2\cos(2\alpha)+|x|^2|y|^2-\langle
x,y\rangle^2 )^2+(|y|^2\sin(2\alpha))^2,\\
\nonumber&&\mathcal {B}:=|y|^2\cos(2\alpha)+|x|^2|y|^2-\langle x,y\rangle^2,\\
\nonumber&&\mathcal {C}:=\langle x,y\rangle \sin(2\alpha), \ \ \  \mathcal{D}:=|x|^4+2|x|^2\cos(2\alpha)+1,
\end{eqnarray}
$0<\alpha<\pi/2$ and $\mathcal{I}m [ \cdot ]$ denote the imaginary part of a complex number.

Based on Berwald's observation (see Lemma \ref{lemBerwald}),
Z. Shen gave the Taylor extensions at the origin $0 \in R^n$ for $x$-analytic projectively flat metrics $F = F(x, y)$ with constant flag
curvature.  He constructed such metrics nearby the origin in $R^n$ using algebraic equations for any given data $F|_{x=0}=\psi(y)$ and $F_{x^k} y^k/(2F)|_{x=0} =\varphi(y)$ \cite{Sh1}.  It is natural to ask if any projectively flat Finsler metric with constant flag curvature is determined by its value at the origin?

We give the positive answer in this paper.
By solving equation (\ref{main}) in real and complex case, we give the classification when $K =0$, $K = -1$ and $K =1$ respectively.
   When $K=0$, we obtain the following.
\begin{thm}\label{thmK=0}
Let $F = F(x,y)$ is a Finsler metric on an open neighborhood $\mathcal{U}$ of the origin in $R^n$. Then $F$
is  projectively flat with zero flag curvature if and only if there exists a Minkowski norm $\psi = \psi(y)$ and
a positively homogeneous function $\phi=\phi(y)$ of degree one on $R^n$ and $C^\infty$ on $R^n \setminus \{0\}$ such that
\begin{equation}\label{thmK=1Eq}
F = \psi ( y + x P) \{ 1 +  P_{y^k} x^k \},
\end{equation}
where  $P = P(x,y)$  satisfies
$ P = \phi (y + x P). $
In this case,
$\psi(y) = F(0,y)$ and $\phi(y) = P(0,y) = F_{x^k} y^k/(2F)|_{x=0}$.
\end{thm}
   Actually, the sufficiency of above theorem is obtained by Z. Shen  in \cite{Sh1}. We prove the necessity in Section 4. It tells us that any
locally projectively flat Finsler metric $F = F(x,y)$ with zero flag curvature is determined by its values ($F(0,y)$ and $P(0,y) = F_{x^k} y^k/(2F)|_{x=0}$) at the origin. The simplest case is that the Euclidean metric $|y|$ can be obtained  by setting $\psi(y) = |y|$ and $\phi(y)=0$. The Berwald's metric (\ref{Berwaldmetric}) can be obtained by setting $\psi(y)=\phi(y)=|y|$.
Actually,
 one can construct many more projectively flat Finsler metrics with zero flag curvature by choosing different $\psi$ and suitable $\phi$.

   The construction of locally projectively flat Finsler metrics  when ${\bf K}=-1$ is somewhat different though they are also determined
   by their behaviors at the origin.  In fact, based on Z. Shen's result (Theorem 1.2 in \cite{Sh1}) and Theorem \ref{mainreal} we prove the following.
\begin{thm}\label{thmK=-1}
Let $F = F(x,y)$ be a Finsler metric on an open neighborhood $\mathcal{U}$ of the origin in $R^n$. Then $F$
is  projectively flat with flag curvature ${\bf K}= - 1$ if and only if  there exists a  Minkowski norm $\psi=\psi(y)$  and a
positively homogeneous function $\phi =\phi(y)$ of degree one on $R^n$ and $C^\infty$ on $R^n \setminus \{0\}$
 such that
\begin{equation} \label{thmK=-1eq}
F(x,y) = \frac{1}{2}\big\{ \Phi_{+} - \Phi_{-} \big\},
\end{equation}
where
\begin{equation} \label{Phiplusminus}
\Phi_{\pm} = (\phi \pm \psi)(y + x \Phi_{\pm}).
\end{equation}
In this case, $\psi(y)=F(0,y)$ and $ \phi(y)= P(0,y) = F_{x^k} y^k/(2F)|_{x=0}$.
\end{thm}
   From this theorem, we can explain why  the known projectively flat Finsler metrics with ${\bf K}=-1$ which can be expressed in elementary functions
are so limited. Actually, it is not easy to solve $\Phi_{\pm}$ in (\ref{Phiplusminus}) with arbitrary $\psi$ and $\phi$. An efficient way is to set
special $\psi$ and $\phi$ such that  (\ref{Phiplusminus}) becomes into a quadratic equation.
For example, by setting $\psi = |y|$ and $\phi = 0$, we get the Riemannian metric $F_{-1}$ with constant section curvature ${\bf K}=-1$. By setting $\psi = \phi = |y|$, we get
\[ F = 2 \frac{ \sqrt{(1-4|x|^2)|y|^2 + 4 \langle x, y \rangle^2}}{1- 4 |x|^2} + 4 \frac{\langle x, y \rangle}{1-4 |x|^2}. \]
By a constant scaling such that $x = \frac{1}{2}x$, we have $\frac{1}{2} F$ is a Funk metric.
More examples are given in \cite{Sh1} in this way.

   In the case when ${\bf K}=1$, we need to express the metrics by the imaginary parts of some complex functions. To solve equation (\ref{main}) in complex case, we need the metric function $F(x,y)$ can be extended to a complex function $F(x,y+x z)$, $z \in \mathbb{C}^n$. Though this excludes some cases, we still have many functions satisfy this condition such as all the analytic functions. And analytic Finsler metric functions  are workable and can be studied directly.
\begin{thm}\label{thmK=1}
Let $F = F(x,y)$ be a Finsler metric on an open neighborhood $\mathcal{U}$ of the origin in $R^n$. Suppose that
$F(x,y)$ can be extended to $F(x,y + x z)$, $z \in \mathbb{C}^n$.
Then
$F$ is projectively flat with flag curvature $\mathbf{K} = 1$ if and only if on $R^n$ there is  a  Minkowski norm $\psi=\psi(y)$  and a
positively homogeneous function $\phi =\phi(y)$ of degree one on $R^n$ and $C^{\infty}$ on $R^n \setminus \{0\}$, and $\psi$ and $\phi$
 can be extended to $\psi(y+x z)$ and $\phi(y+x z)$ ($z \in \mathbb{C}^n$)
such that
\begin{equation} \label{HF}
F(x, y) = \mathcal{I}m [ \Psi(x,y) ],
\end{equation}
where
\begin{equation}\label{H}
\Psi = \phi (y + x \Psi ) + i \psi(y + x \Psi).
\end{equation}
In this case,  $\psi(y)=F(0,y)$ and $ \phi(y)= P(0,y) = F_{x^k} y^k/(2F)|_{x=0}$.
\end{thm}
It is not easy to give the expression of $\Psi$ in (\ref{H}) for most of the choices of $\phi$ and $\psi$.
However, we can get some special ones by choosing suitable $\psi$ and $\phi$.
For example, the Riemannian metric $F_{+1}$ with constant sectional curvature ${\bf K}=1$ can be obtained
by setting $\psi=|y|$ and $\phi=0$.
By setting $\phi+ i\psi=i e^{-i\alpha}|y|$, Z. Shen obtained
Bryant's metric (\ref{Bryant}). Recently, by setting $\phi=\phi=|y|+\langle a, y\rangle$ (a is a constant vector), we get
another metric with double square roots which is projectively flat and with flag curvature $K=1$ \cite{XuLi}.
In this paper, we
construct a new locally projectively flat Finsler metric with constant flag curvature ${\bf K}=1$. Its
 $\psi$ and $\phi$  both are with double square roots. See Example \ref{newex}.

   From Theorem \ref{thmK=0} - \ref{thmK=1}, we can see that locally projectively flat Finsler metrics with
 constant flag curvature are totally determined by its behaviors at the origin. Any pair of $\psi=\psi(y)$ and $\phi=\phi(y)$
  can produce locally projectively flat Finsler metrics with constant flag curvature ${\bf K} =0$, $-1$ or $+1$ in three
  different ways and vice versa.

   In recent years, many Finsler metrics composed of Riemannian metrics and $1$-forms are studied such as
(general) ($\alpha$, $\beta$)-metrics, spherical symmetric Finsler metrics, and etc.
In 2006, we classified locally projectively flat ($\alpha$, $\beta$)-metrics \cite{LiShen} into tree types.
In 2012, L. Zhou give the classification of projectively flat spherically symmetric  Finsler metrics with constant
flag curvature \cite{Zhou}. His proof based on complicated computation and related analysis on some PDEs.
As an application of Theorem \ref{thmK=0} - \ref{thmK=1},
we give the classification of spherically symmetric Finsler metrics much easier in Section \ref{app}.

\section{Preliminaries}

   A {\it Minkowski norm}  on a vector space is a $C^{\infty}$
 function  $\psi : V \setminus \{ 0 \} \rightarrow [0, +\infty)$ satisfying: (i) $\psi (y) =0$ if and only if $y=0$;
 (ii) $\psi$ is positively homogeneous of degree one, i.e., $\psi(\lambda y) = \lambda \psi$, $\lambda > 0$;
 (iii) $\psi(y)$ is strongly convex, i.e.,  the matrix $g_{ij}(y) := [ \frac{1}{2} F^2]_{y^i y^j} (y)$ is positive definite.
 A {\it Finsler metric} $F = F(x,y)$ on a  manifold $M$ is a $C^{\infty}$ function on $TM \setminus \{ 0 \}$ such
 that $F|_{T_{x}M}$ is a Minkowski  norm on $T_{x}M$ for each $x\in M$.

   Consider a Finsler metric $F=F(x,y)$ on an open domain $\mathcal{U} \subset R^n$. The geodesics of $F$ are determined by the following ODEs:
\[ \ddot{x} + 2 G^i (x, \dot{x}) = 0,\]
where $ G^i = G^i(x,y)$ are called {\it geodesic coefficients }
given by
\[ G^i = \frac{1}{4} g^{il} \Big\{  [F^2]_{x^m y^l} y^m -[F^2]_{x^l} \Big\}. \]

  As an extension of sectional curvature in Riemann geometry, for each tangent plane $\Pi\subset T_xM$ and $y\in \Pi$, the
\emph{flag curvature} of $(\Pi,y)$ is defined by
\begin{eqnarray}
\nonumber
{\bf K}(\Pi,y)=\frac{g_{im}R^i_{\ k} u^ku^m}{F^2g_{ij}u^iu^j-[g_{ij}y^iu^j]^2},
\end{eqnarray}
where $\Pi=span\{y,u\}$, and
\begin{eqnarray}
\nonumber R^i_{\ k}=2\frac{\partial G^i}{\partial
x^k}-y^j\frac{\partial^2G^i}{\partial x^j\partial
y^k}+2G^j\frac{\partial^2G^i}{\partial y^j\partial
y^k}-\frac{\partial G^i}{\partial y^j}\frac{\partial G^j}{\partial
y^k}.
\end{eqnarray}
Finsler metric $F$ is of {\it scalar flag curvature} if its flag
curvature ${\bf K}(\Pi,y)={\bf K}(x,y)$ is independent of tangent plane $\Pi$.
If $F$ is a  Riemannian metric, the flag curvature ${\bf K}(\Pi,y)={\bf K}(\Pi)$ is
independent of $y$. Finsler metric $F$ is said to be with {\it
constant flag curvature} if ${\bf K}=\lambda$ is a constant. In this case,
\[ R^i_{\ k}=\lambda\{F^2\delta^i_k-FF_{y^k}y^i\}.\]

  $F$ is said to be \textit{projectively flat} in
$\mathcal{U}$ if all geodesics are straight lines. This is equivalent to $G^i = P(x,y) y^i$, where $P =F_{x^k}y^k / (2F)$ is called the {\it projective factor} of $F$. In 1903, G. Hamel proved that $F$ is locally projectively flat if and only if
\begin{equation}\label{proj}
F_{x^k} - F_{x^l y^k} y^l = 0.
\end{equation}
In this case, the flag curvature of $F$  is a scalar function on $T{\cal U}$ given by
\begin{equation}
{\bf K} = \frac{P^2 - P_{x^m} y^m }{F^2}.\label{flagC}
\end{equation}
This observation is  due to L. Berwald \cite{Be2}. In his paper, he proved the following lemma.
\begin{lem}\label{lemBerwald}
Let $F=F(x,y)$ be a Finsler metric on an open subset $\mathcal{U} \subset R^n$. Then $F$ is projectively flat if and only if
there is a positively $y$-homogeneous function of degree one, $P = P(x, y)$, and a
positively homogeneous function of degree zero, $\mathbf{K} = \mathbf{K}(x, y)$, on $T \mathcal{U}  \simeq \mathcal{U} \times R^n$
such that
\begin{equation}\label{Berwald1}
F_{x^k} = (P F)_{y^k},
\end{equation}
\begin{equation} \label{Berwald2}
P_{x^k} = P P_{y^k} - \frac{1}{3 F} (\mathbf{K} F^3)_{y^k}.
\end{equation}
In this case, $P$ is the projective factor of $F$.
\end{lem}
It is easy to see that if ${\bf K}=0$ then the projective factor $P$ satisfies
\begin{equation} \label{Pequation}
P_{x^k} = P P_{y^k}.
\end{equation}

In the case $\mathbf{K} = \lambda \neq 0$, L. Berwald discovered the following lemma.
\begin{lem}\label{lemBerwald1}
Let $F=F(x,y)$ be a Finsler metric on an open subset $\mathcal{U} \subset R^n$. Then $F$ is projectively flat with constant flag
curvature $\mathbf{K} = \lambda \neq 0$ if and only if
\begin{equation}\label{mainequation}
(\Phi_{\pm})_{x^k} = \Phi_{\pm} (\Phi_{\pm})_{y^k}.
\end{equation}
where $P = P(x,y)$ is the projective factor of $F$ and
\[ \Phi_{\pm} = P \pm \sqrt{-\lambda} F. \]
\end{lem}
Thus the key problem to classify locally projectively flat Finsler metrics with constant flag curvature is  to solve equation (\ref{Pequation}) and (\ref{mainequation}).

\section{Solution of $\Phi_{x^k} = \Phi \Phi_{y^k}$}
It is difficult to solve  (\ref{Pequation}) and (\ref{mainequation}) directly for their nonlinearity.
In 2003, inspired by the structure of Funk metric, Z. Shen find a solution of (\ref{Pequation}) as  following.
\begin{lem}\label{lemShen}(\cite{Sh1}) Let $\phi = \phi(y)$ be an arbitrary  positively homogeneous function of degree one one $R^n$. Suppose that
$\phi$ is $C^\infty$ on $R^n \setminus \{0\}$. Then there is a unique real-valued
function $\Phi = \Phi(x, y)$ satisfying the following
\[ \Phi = \phi( y + x \Phi). \]
Moreover, $\Phi$ satisfies
 \[ \Phi_{x^k} = \Phi \Phi_{y^k}. \]
\end{lem}
Then a natural problem is to determine all the solutions. Is there any other solution?
We prove the following lemma and show that there is no other solutions.
\begin{lem}\label{lemmainreal}
Let $\Phi = \Phi(x,y)$ be a positively $y$-homogeneous  function of degree one. Suppose $\Phi$ is $C^\infty$ on $T\mathcal{ U }\setminus \{0\}= \mathcal{U} \times R^n \setminus \{0\}$ satisfying $\Phi_{x^k} = \Phi \Phi_{y^k}$, where $\mathcal{U}$ is an open neighborhood of the origin in $R^n$. Then
there is a unique positively homogeneous function $\phi = \phi(y)$ of degree one on $R^n$ such that
\begin{equation} \label{Pphi}
\Phi = \phi(y + x \Phi).
\end{equation}
In this case, $\phi(y)= \Phi(0,y)$.
\end{lem}
{\bf Proof:} Let $f(t) := t - \Phi(x, y - x t)$.
Fixing $x$ and $y$, we need to prove that locally there is a unique $t_{o}$ such that
$f(t_{o})=0$. We divide the proof in two cases.

{\bf Case (i)} $y - x t \neq 0$ for any $t$.

Observe that there is a small $\epsilon  > 0 $ such that
for any $x \in R^n$ with $|x| < \epsilon$, at any $t$ where $y - x t \neq 0$,
\[ | \Phi_{y^k}x^k | \leq \sup_{|\eta|=1} \Phi_{y^k}(x,\eta)|x| \leq \epsilon \sup_{|\eta|=1} \Phi_{y^k}(x,\eta) \]
\be \label{mainrealine1}
 \frac{3}{2}\geq f'(t) = 1 + \Phi_{y^k}(x, y-xt)x^k \geq \frac{1}{2}.
 \ee

By mean value theorem, for any $\bar{t}$, there is a $\xi \in (t, \bar{t})$ (or $\xi \in (\bar{t}, t)$) such that
\[ f(t) - f(\bar{t}) = f'(\xi) (t-\bar{t}). \]
Then by above equation and (\ref{mainrealine1}), we get
$f(t)$ is a monotonic increasing function satisfying $f(t)\rightarrow + \infty (t\rightarrow +\infty)$ and $f(t)\rightarrow - \infty (t\rightarrow -\infty)$.
Thus there is a unique $t_{o}$ such that $f(t_{o}) = 0$.

{\bf Case (ii) }$y - x \lambda =0$ for some $\lambda$.

In this case,
\[ f(t) = t - \Phi(x, (\lambda -t)x ) = t - |\lambda - t| \Phi(x,x).\]
Then there is a small $\epsilon  > 0 $ such that
for any $x \in R^n$ with $|x| < \epsilon$, $f(t)$ is a monotonic increasing function satisfying $f(t)\rightarrow + \infty (t\rightarrow +\infty)$ and $f(t)\rightarrow - \infty (t\rightarrow -\infty)$.
Thus there is a unique $t_{o}$ such that $f(t_{o}) = 0$.

Then we get the unique solution by setting $\phi(x,y) = t_{o}$ such that
\begin{equation} \label{phiPhi}
\phi(x,y) = \Phi(x, y - x \phi(x,y)).
\end{equation}
Next we prove $\phi(x,y)$ is independent of $x$.
Set
\begin{equation}\label{Ppsi}
\eta = y - x \phi.
\end{equation}
Then differentiating (\ref{phiPhi}) with respect to $y^k$ and $x^k$ respectively,
we get
\begin{equation}\label{psiy}
(1+ \Phi_{\eta^l} x^l)\phi_{y^k} = \Phi_{\eta^k},
\end{equation}
\begin{equation}\label{psix}
(1+ \Phi_{\eta^l} x^l)\phi_{x^k} = \Phi_{x^k} - \Phi_{\eta^k} \phi = \Phi_{x^k} - \Phi_{\eta^k} \Phi =0.
\end{equation}
Here the assumption $\Phi_{x^k} = \Phi \Phi_{y^k}$ is used.
If $\Phi_{\eta^l} x^l = -1 $, then by (\ref{psiy}) $\Phi_{\eta^k} =0$. It is a contradiction. Then by (\ref{psix}), we obtain
  $\phi=\phi(y)$.
\qed

Then by Lemma \ref{lemShen} and Lemma \ref{lemmainreal} we obtain the following.
\begin{thm}\label{mainreal}
Let $\Phi = \Phi(x,y)$ be a positively $y$-homogeneous function of degree one on $T\mathcal{ U }= \mathcal{U} \times R^n$, where $\mathcal{U}$ is an open neighborhood of the origin in $R^n$. Suppose that $\Phi$ is $C^{\infty}$ on $T\mathcal{ U }\setminus \{0\}= \mathcal{U} \times R^n \setminus \{0\}$.
Then
\be \label{maineqlem}
 \Phi_{x^k} = \Phi \Phi_k \ee
 if and only if there is a unique positively homogeneous function $\phi = \phi(y)$ of degree one on $R^n$ and
 $C^{\infty}$ on $R^n \setminus \{0\}$ such that
\begin{equation} \label{Pphi}
\Phi = \phi(y + x \Phi).
\end{equation}
In this case, $\phi(y)=\Phi(0,y)$.
\end{thm}
This theorem tell us that each solution of (\ref{maineqlem}) corresponding to a unique positively homogeneous function.
It plays an important role in our proofs of Theorem \ref{thmK=0}, \ref{thmK=-1}.

\section{{\bf K}=0}

In this section, we are going to determine the structure of projectively flat Finsler metrics with zero flag curvature. In \cite{Sh1}, Z. Shen construct some
examples based on following theorem.
\begin{thm}\label{thmShenK=0}(\cite{Sh1})
Let $\psi(y)$ be an arbitrary Minkowski norm on $R^n$ and $\phi(y)$ be an
arbitrary positively homogeneous function of degree one on $R^n$. Define $P(x, y)$
by
\[ P(x, y) = \phi(y + x P(x, y)). \]
Let
\[ F = \psi ( y + x P) \{ 1 +  P_{y^k} x^k \}.\]
Then $F(x,y)$ is a locally projectively flat Finsler metric with zero flag curvature.
\end{thm}
In fact, all known locally projectively flat Finsler metrics with zero flag curvature can be determined by this theorem.
It leads us to study whether all such metrics are determined in this way or not. By Lemma \ref{lemBerwald}, (\ref{Pequation}) and
Theorem \ref{mainreal}, we have that the projective factor of  any locally projectively flat Finsler metric $F=F(x,y)$ with ${\bf K}=0$ must be uniquely determined
by a positively homogeneous function of degree one on $R^n$. Then by Theorem \ref{thmShenK=0}, we can construct a projectively flat Finsler metric
$\tilde{F}=\tilde{F}(x,y)$ with ${\bf K}=0$ whose projective factor are same with the one of $F$. To tell the relation between these two metrics, we prove the following lemma.
\begin{lem} \label{samefactor}
Let $F=F(x,y)$ and $\tilde{F} = \tilde{F}(x,y)$ are two locally projectively flat Finsler metrics. If they have the same projective factor $P = P(x,y)$,
then one of the following holds:

i)  $\tilde{F} = c F(x,y)$, where $c$ is a positive constant;

ii)  $F$ and $\tilde{F}$ both have zero flag curvature and $\tilde{F} = \Theta F(x,y)$, where $\Theta=\Theta(x,y)$ satisfies
\begin{equation} \label{ThetaPrelation}
\Theta_{x^k} = P \Theta_{y^k}.
\end{equation}
\end{lem}
{\bf Proof:} By the definition of projective factor and the assumption, we have
\[ \frac{F_{x^k} y^k}{F} = \frac{\tilde{F}_{x^k} y^k}{\tilde{F}} = 2 P. \]
Let $\tilde{F} = \Theta F$. Then by above equation we get
\begin{equation}\label{Theta0}
\Theta_{x^k}y^k = 0.
\end{equation}
Differentiating it respect to $y^k$ yields
\begin{equation}\label{Theta1}
\Theta_{x^k} + \Theta_{x^l y^k} y^l = 0.
\end{equation}
By assumption $F$ and $\tilde{F} = \Theta F$ both are locally projectively flat,
substituting $F$ and $\Theta F$ into G. Hamel's equation (\ref{proj}) yields
\begin{equation} \label{GF}
F_{x^k} - F_{y^k x^l}y^l =0,
\end{equation}
\begin{equation} \label{GThetaF}
\Theta ( F_{x^k} - F_{y^k x^l}y^l ) + F (\Theta_{x^k} - \Theta_{y^k x^l}y^l) -\Theta_{y^k} F_{x^l}y^l -  F_{y^k} \Theta_{x^l}y^l =0.
\end{equation}
Substituting (\ref{Theta0}), (\ref{Theta1}) and (\ref{GF}) into (\ref{GThetaF}) yields

\[ 2 \Theta_{x^k} F - \Theta_{y^k} F_{x^l} y^l =0. \]
Then by the definition of projective factor we have
\begin{equation} \label{ThetaP1}
\Theta_{x^k} = P \Theta_{y^k}.
\end{equation}
Differentiating above equation with respective to $x^l$ and contracting with $y^l$ yields
\[  \Theta_{x^k x^l} y^l = P_{x^l}y^l \Theta_{y^k} + P \Theta_{y^k x^l}y^l.\]
By (\ref{Theta0}) we have $ \Theta_{x^k x^l} y^l =0$. Substituting it and (\ref{Theta1}) into above equation we get
\[ P \Theta_{x^k}  = P_{x^l}y^l \Theta_{y^k}.\]
Then by (\ref{ThetaP1}) we get
\[ \Theta = constant \ \ \  \text{or} \ \ P_{x^k}y^k = P^2.  \]
In the latter case
by (\ref{Berwald2}) in Lemma \ref{lemBerwald}, we get the flag curvatures of $F$ and $\tilde{F}$ both are zero.
\qed

To prove Lemma \ref{SolutionTheta} for solving (\ref{ThetaPrelation}), we need
\begin{lem}\label{lemimplicit}
Let
\begin{equation}
E^k := \xi^k - y^k - x^k P(x,y) =0,
\end{equation}
where $P = P(x,y)$ is a positively $y$-homogeneous  function of degree zero on $T\mathcal{ U }= \mathcal{U} \times R^n$ and $C^\infty$ on $T\mathcal{ U } \setminus \{0\}= \mathcal{U} \times R^n \setminus \{0\}$ satisfying
\begin{equation} \label{partialypartialx}
P_{x^k} = P P_{y^k}.
\end{equation}
 Then
\be \label{implicit}
 \frac{\partial y^k}{\partial x^j} = - P \delta^{k}_{j}
 \ee
\end{lem}
{\bf Proof:} By a direct computation, we have
\[ E^k_{y^j} = - \delta^k_j - x^k P_{y^j}, \]
\begin{equation}
E^k_{x^j} = - P \delta^k_j - x^k P_{x^j}
          = - P (\delta^k_j + x^k P_{y^j}).
\end{equation}
By implicit differentiation and (\ref{partialypartialx}), we get (\ref{implicit}).
\qed

\begin{lem}\label{SolutionTheta}
Let $\Theta=\Theta(x,y)$  is a positively $y$-homogeneous  function of degree zero on $T\mathcal{ U } = \mathcal{U} \times R^n $ and $C^{\infty}$ on $T\mathcal{ U } \setminus \{0\}= \mathcal{U} \times R^n \setminus \{0\} $. If $\Theta$ satisfies
\begin{equation} \label{ThetaP}
\Theta_{x^k} = P \Theta_k,
\end{equation}
where $P = P(x,y)$ is a positively $y$-homogeneous function of degree one such that $P_{x^k} = P P_{y^k}$,
then
\[ \Theta = \psi(y + x P), \]
where $\psi = \psi(y)$ is a positively $y$-homogeneous function of degree zero.
\end{lem}
{\bf Proof:} Let
\begin{equation} \label{newLemEq1}
\xi^k = y^k + x^k P(x,y).
\end{equation}
Regarding $y$ as a function of $x$ and $\xi$, then by Lemma \ref{lemimplicit} we have
\begin{equation} \label{newLemEqpartial}
\frac{\partial y^l }{\partial x^k} = - P \delta^l_k.
\end{equation}
Let
\[ \psi(x, \xi) = \Theta(x,y). \]
We only need to prove $\psi(x, \xi) = \psi(\xi)$, i.e., $\psi_{x^k} =0$. In fact,
\begin{equation}
\begin{split}
\psi_{x^k} =& \Theta_{x^k} + \Theta_{y^l}\frac{\partial y^l}{\partial x^k} \\
           =&  \Theta_{x^k} - P \Theta_{y^l}  =0.
\end{split}
\end{equation}
Here we used (\ref{newLemEqpartial}) and (\ref{ThetaP}). Thus
\[ \Theta(x,y) = \psi(\xi) =\psi(y+xP).\]
\qed

{\bf Proof of Theorem \ref{thmK=0}: } The sufficiency is obtained by Z. Shen's Theorem \ref{thmShenK=0}. We only need to prove the necessity.
If $F=F(x,y)$ is a locally projectively flat Finsler metric with ${\bf K}=0$ with its projective factor $P = P(x,y)$, then by Theorem \ref{mainreal} $P$ is uniquely determined by  a positively homogeneous function $\phi = \phi(y)$ of degree one on $R^n$. By Theorem \ref{thmShenK=0}, we can construct a projectively flat Finsler metric $\tilde{F}$ by any Minkowski norm $\bar{\psi}=\bar{\psi}(y)$, i.e.,
\[ \tilde{F} = \bar{\psi}(y + x P)  \{ 1 +  P_{y^k} x^k \}. \] 
Then
by Lemma \ref{SolutionTheta} and Lemma \ref{samefactor}, there is a  positively $y$-homogeneous function $\tilde{\psi}=\tilde{\psi}(y)$ of degree zero such that
\[ F = \tilde{\psi}(y + x P) \tilde{F} = \tilde{\psi}(y + x P) \bar{\psi}(y + x P)  \{ 1 +  P_{y^k} x^k \}. \]
By setting $\psi = \tilde{\psi}\bar{\psi}$ we get
\[ F = \psi(y + x P)  \{ 1 +  P_{y^k} x^k \}. \]
\qed

\section{$ {\bf K} = -1$}
The construction of locally projectively flat Finsler metrics with ${\bf K}=-1$ is different from the case when ${\bf K}=0$.
By (\ref{mainequation}), we have
\begin{equation}\label{K=-1(1)}
(P + F)_{x^k} = (P + F) (P + F)_{y^k},
\end{equation}
\begin{equation}\label{K=-1(2)}
(P - F)_{x^k} = (P - F) (P - F)_{y^k}.
\end{equation}
In 2003, Z. Shen constructs some metrics based on Theorem 1.2 in \cite{Sh1}. Actually Theorem 1.2 in \cite{Sh1} is the sufficiency of
our Theorem \ref{thmK=-1}.

{\bf Proof of Theorem \ref{thmK=-1}:} We only need to proof the necessity.  By (\ref{K=-1(1)}) and (\ref{K=-1(2)}), $P + F$ and $P-F$ satisfy equation (\ref{Pphi}) in Theorem \ref{mainreal} respectively. Then by Theorem \ref{mainreal} there exist unique $\phi_{+} =\phi_{+}(y)= (P+F)|_{x=0}$ and $\phi_{-}=\phi_{-}(y)=(P-F)|_{x=0}$ such that
\begin{equation} \label{thmK=-1_PplusF}
P + F = \phi_{+} (y + x (P+F)),
\end{equation}
\begin{equation} \label{thmK=-1_PminusF}
P - F = \phi_{-} (y + x (P-F)).
\end{equation}
Setting $F(0,y) = \psi(y)$ and $P(0,y)=\varphi(y)$, we have
\[ \phi_{+} = \varphi + \psi, \ \ \ \ \phi_{-} = \varphi - \psi.\]
Then (\ref{Phiplusminus}) is just (\ref{thmK=-1_PplusF}) and (\ref{thmK=-1_PminusF}). Thus (\ref{thmK=-1eq}) is obtained.
\qed

\section{$ K = 1$}

In this case, (\ref{mainequation}) is equivalent to
\begin{equation}\label{eqcomplex}
( P + i F )_{x^k} = ( P + i F ) (P + i F)_{y^k}.
\end{equation}
If we still want to express the solutions of above equation in an algebraic way, we need to add some condition
on $P=P(x,y)$ and $F(x,y)$ such that they can be extended to $\mathcal{U} \times C^n$, $\mathcal{U} \subset R^n$.
It is easy to
see that if they are $y$-analytic, then they can be extended.

By the similar argument as in Lemma \ref{lemmainreal}, we have
\begin{thm} \label{lemcomplex}
Let $\Psi = P + i F$, where $P = P(x,y)$ and $F = F(x,y)$ are two positively $y$-homogeneous functions of degree one on $T\mathcal{ U }= \mathcal{U} \times R^n $. Suppose $P$ and $F$ are both $C^\infty$ on $T\mathcal{ U }\setminus \{0\}= \mathcal{U} \times R^n \setminus \{0\}$ and  can be extended to $\mathcal{U} \times \mathbb{C}^n$.
Then
\be \label{complex0}
\Psi_{x^k} = \Psi\Psi_{y^k}
\ee
 if and only if
there are two  positively homogeneous functions $\phi = \phi(y)$ and $\psi(y)$ of degree one on $R^n$ and $C^\infty$ on $R^n \setminus \{0\}$ which can be extended to $C^n$
such that
\begin{equation} \label{complex}
\Psi = \phi(y + x \Psi) + i \psi(y+ x \Psi).
\end{equation}
\end{thm}
{\bf Proof:} The sufficiency is first discovered by Z. Shen. It can be verified directly in the same way in Lemma  \ref{lemShen}.
We only need to prove the necessity. The proof is similar as in Lemma  \ref{mainreal}. We give the main part here to prove the following functions
$f(t)$ and $g(t)$ have unique zero point.
Set real function
\[ f(t):= t - \mathcal{R}e\Big[ P(x, y-x(t+i s)) + i F(x, y-x(t+is)) \Big]. \]
Then there is a small $\epsilon_1  > 0 $ such that
for any $x \in R^n$ with $|x| < \epsilon_1$, at any $t$ where $y - x (t+is) \neq 0$,
\[ f'(t) = 1 + \mathcal{R}e\Big[ P_{y^k}(x, y-x(t+i s))x^k + i F_{y^k}(x, y-x(t+is))x^k \Big] \geq \frac{1}{2}. \]
Thus there is a unique $t_{o}=t_{o}(s)$ such that $f(t_{o}(s))=0.$

Similarly, set real function
\[ g(s):= s - \mathcal{I}m\Big[ P(x, y-x(t_{o}(s)+i s)) + i F(x, y-x(t_{o}(s)+is)) \Big]. \]
Then there is a small $\epsilon_2  > 0 $ such that
for any $x \in R^n$ with $|x| < \epsilon_2$, at any $s$ where $y - x (t_{o}(s)+is) \neq 0$,
\be
\begin{split}
g'(s) &= 1 + \mathcal{I}m\Big[ x^k P_{y^k}(x, y-xt_{o}(s) - i x s)(t_{o}'(s)+i) + i x^k F_{y^k}(x, y-xt_{o}(s) -i x s)(t_{o}'(s)+i) \Big] \\
& \geq \frac{1}{2}.
\end{split}
\ee
Thus there is a unique $s_{o}$ such that $g(s_{o})=0.$ Then we get the unique solution by setting
$\phi(x,y)=t_{o}(s_{o})$ and $\psi(x,y)=s_{o}$ such that
\be \label{clemEq}
\begin{split}
\phi(x,y) + i\psi(x,y) &= P(x, y-x\phi(x,y)-ix\psi(x,y))+i F(x, y-x\phi(x,y)-ix\psi(x,y))\\
&=\Psi(x,y-x\phi(x,y)-ix\psi(x,y)).
\end{split}
\ee
To prove $\phi(x,y)+ i\psi(x,y)$ is independent of $x$, we set
\[ \eta =  y - x \phi -i x \psi.\]
Differentiating (\ref{clemEq}) with respect to $y^k$ and $x^k$ respectively yields
\be \label{clem1}
(1+\Psi_{\eta^l}x^l)(\phi_{y^k}+i \psi_{y^k}) = \Psi_{\eta^k},
\ee
\be \label{clem2}
(1+\Psi_{\eta^l}x^l)(\phi_{x^k}+i \psi_{x^k}) = \Psi_{x^k} - \Psi_{\eta^k}(\phi + i \psi) = \Psi_{x^k} - \Psi_{\eta^k}\Psi =0.
\ee
Here the assumption  $\Psi_{x^k} = \Psi\Psi_{\eta^k}$ is used in (\ref{clem2}). If $\Psi_{\eta^l} x^l = -1 $, then by (\ref{clem1}) $\Psi_{\eta^k} =0$.
It is a contradiction. Then by (\ref{clem2}), we obtain
  \[ \psi=\psi(y), \ \ \ \phi=\phi(y). \]
\qed

By above Theorem and (\ref{eqcomplex}) we give the proof of Theorem \ref{thmK=1}.

{\bf Proof of Theorem \ref{thmK=1}}: The sufficiency was first discussed in \cite{Sh1} and can be verified directly.
 We only need to proof the necessity. By (\ref{eqcomplex}), we have that $\Psi = P+ i F$ satisfies (\ref{complex0}) in Theorem \ref{lemcomplex}.
 Then by Thereom \ref{lemcomplex} there exist two
positively homogeneous functions $\psi =\psi(y)$ and $\phi=\phi(y)$ of degree one on $R^n$ such that
$\Psi$ satisfies (\ref{H}). In this case,
$F(x,y) = \mathcal{I}m [\Psi(x,y)]$ and
$\psi(y) = F(0,y)$ must be a Minkowskian norm.
\qed

Based on Theorem \ref{thmK=1}, we can construct a new projective flat Finsler metric with constant flag curvature ${\bf K}=1$ from a Minkowskia norm composed of double square roots.
\begin{ex}\label{newex}
Let $\mathcal{U}$ is an open neighborhood at the origin in $R^n$ and  $\mathcal{\tilde{U}}$ is an open neighborhood at the origin in $R^m$. Set
\[ \phi = \frac{\sqrt{2}}{2} \sqrt{ \sqrt{|y|^4 + |\tilde{y}|^4 } - |y|^2 },  \]
\[ \psi = \frac{\sqrt{2}}{2}\sqrt{ \sqrt{|y|^4 + |\tilde{y}|^4 } + |y|^2 },  \]
where $|y|$ and $|\tilde{y}|$ are Euclidean norms on $\mathcal{U}$ and $\mathcal{\tilde{U}}$ respectively. It can be verified directly that
$\psi$ is a Minkowski norm.
Then by (\ref{HF}) in Theorem \ref{thmK=1}
\be \label{newexample}
 F = \mathcal{I}m \Big[ \frac{- \langle x,y\rangle+i \langle \tilde{x},\tilde{y}\rangle  + i \sqrt{ (|y|^2-i|\tilde{y}|^2)(1+|x|^2-i|\tilde{x}|^2)-(\langle x,y\rangle- i\langle \tilde{x},\tilde{y}\rangle)^2 } }{1+|x|^2-i|\tilde{x}|^2} \Big]
\ee
is a projectively flat Finsler metric with constant flag curvature ${\bf K} =1$.
It is easy to see that on $\mathcal{U}$ it is the Riemannian metric $F_{+1}$ in (\ref{Rie}).

\end{ex}

\section{Applications} \label{app}
   Let $\Omega\subseteq R^n$ be a convex domain. A Finsler metric $F = F(x,y)$ on $\Omega$ is called a {\it spherically symmetric} Finsler metric if
   \[ F(S x, S y) = F(x,y), \]
   for all $S \in O(n)$. Obviously, many known special Finsler metrics are spherically symmetric metrics such as
   (\ref{Rie}), (\ref{Funk}), (\ref{Berwaldmetric}) and  (\ref{Bryant}). It is proved  that  any spherically symmetric Finsler metric $F$ can be expressed by
   \[ F(x,y) = |y| \zeta(|x|, \frac{\langle x, y\rangle}{|y|}),\]
   where $\zeta = \zeta(s,t)$ is a $C^{\infty}$ function  \cite{Mo}.
In \cite{Zhou}, L. Zhou studies projectively flat spherically symmetric Finsler metrics with constant flag curvature and give the classification by long computation
and some analysis on related PDEs. Now, by Theorem \ref{thmK=0}, \ref{thmK=-1} and \ref{thmK=1} we can give the classification much easier.

\begin{lem} \label{lemsym}
Let $F = |y|\zeta(|x|,\frac{\langle x,y\rangle}{|y|})$ be a projectively flat spherically symmetric Finsler metric on a convex domain $\Omega\subseteq R^n$.
Then
by a constant scaling
\[ F(0,y) =|y| \]
and
\[ P(0,y) = \frac{F_{x^k}y^k}{2F}|_{x=0}= c|y|, \]
where $c$ is a constant.
\end{lem}
{\bf Proof:} It is obvious by the definition of spherically symmetric Finsler metric.

\begin{thm}
Let $F = |y| \zeta(|x|,\frac{\langle x,y\rangle}{|y|})$ be a spherically symmetric Finsler metric on on a convex domain $\Omega\subseteq R^n$. Then $F$ is locally projectively flat with zero flag curvature
 if and only if

(i) $F = |y|$; or

(ii) \be \label{symK=0F}
F = \frac{|y|^4}{z (\langle x,y\rangle \pm z)^2}. \ee
where $z = \sqrt{(1-c^2|x|^2)|y|^2+c^2 \langle x,y \rangle^2}$, $c$ is a nonzero constant.
\end{thm}
{\bf Proof:} The sufficiency can be verified directly. We only need to prove the necessity.
By Theorem \ref{thmK=0}, a key problem is to determin $F(0,y)$ and $P(0,y)= F_{x^k} y^k/(2F)|_{x=0}$.
By Lemma \ref{lemsym}, by a constant scaling on $y$
\[ F(0,y) =  |y|. \]
 By Theorem \ref{thmK=0}, to get the projective factor $P=P(x,y)$, we only need to solve the equation
\[ P = \varphi (y + x P), \]
where $\varphi = \varphi(y)$ is an arbitrary positively homogeneous function of degree one on $y$. Obviously, the only positively $y$-homogeneous function of degree one in this case is  $c |y|$, where  $c$  is a  constant.
Then
\be \label{ssK=0P}
 P = c | y + x P |.
\ee
If $c=0$, then $P =0$.
If $c\neq 0$, then by solving (\ref{ssK=0P})
we get
\[ P = \frac{c^2<x,y> + sgn(c) \sqrt{ c^2(1-c^2 |x|^2)|y|^2+ c^4 <x,y>^2}}{1-c^2 |x|^2}. \]
Then by (\ref{thmK=1Eq}), we obtain (\ref{symK=0F}).
\qed

When ${\bf K}=-1$, the proof is similar. We just need to use Theorem \ref{thmK=-1} here.
\begin{thm}
Let $F = |y| \zeta(|x|,\frac{\langle x,y\rangle}{|y|})$ be a spherically symmetric Finsler metric on a convex domain $\Omega\subseteq R^n$. Then $F$ is locally projectively flat with constant flag curvature
${\bf K}=-1$ if and only if
\be \label{symK=-1F}
\begin{split}
F =& \frac{1}{2}\Big\{ \frac{(c+1)^2\langle x,y\rangle + sgn(c+1) \sqrt{ (c+1)^2(1-(c+1)^2 |x|^2)|y|^2+ (c+1)^4 \langle x,y\rangle^2}}{1-(c+1)^2 |x|^2} \\
&- \frac{(c-1)^2\langle x,y\rangle + sgn(c-1) \sqrt{ (c-1)^2(1-(c-1)^2 |x|^2)|y|^2+ (c-1)^4 \langle x,y\rangle^2}}{1-(c-1)^2 |x|^2}\Big\},
\end{split}
\ee
where $c$ is a constant.
\end{thm}
{\bf Proof:} The sufficiency can be verified directly. We only need to prove the necessity. By Lemma \ref{lemsym},
\[ F(0,y) =  |y|,\ \ \  P(0,y) = c |y|, \]
where $c$ is a constant.
Then
\be
 \Phi_{+} =  (c+1)|y + x \Phi_{+}|,
\ee
\be
 \Phi_{-} =  (c-1)|y + x \Phi_{-}|.
\ee
Solving above two equations, we get
\[ \Phi_{+} =  \frac{(c+1)^2\langle x,y\rangle + sgn(c+1) \sqrt{ (c+1)^2(1-(c+1)^2 |x|^2)|y|^2+ (c+1)^4 \langle x,y\rangle^2}}{1-(c+1)^2 |x|^2},  \]
\[ \Phi_{-} =  \frac{(c-1)^2\langle x,y\rangle + sgn(c-1) \sqrt{ (c-1)^2(1-(c-1)^2 |x|^2)|y|^2+ (c-1)^4 \langle x,y\rangle^2}}{1-(c-1)^2 |x|^2}.  \]
Then by (\ref{thmK=-1eq}) we get (\ref{symK=-1F}).
\qed

In \cite{Zhou}, L. Zhou claimed a "new" projectively flat Finsler metric (\ref{Zhoumetric}) with ${\bf K} = -1$ is found.
Actually, we can prove it is also can be written as (\ref{symK=-1F}).

\begin{ex}
Let $\Omega = B^n (\sqrt{2(d_2 - d_1)}) \subset R^n$  with Finsler metric
\be \label{Zhoumetric}
 F = \frac{|y| c_1(z_1)}{ c_1(z_1)^2 - (z_2 + c_2(z_1)^2)^2},
 \ee
where
\[ z_1 := \sqrt{|x|^2 - \frac{<x,y>^2}{|y|^2}}, \ \ \ \ z_2 := \frac{\langle x,y\rangle}{|y|}, \]
\[ c_1(z_1) := \frac{\sqrt{2}}{2}\sqrt{2d_2 - z_1^2 + \sqrt{(2d_2 - z_1^2)^2 - 4d_1^2}}, \]
\[ c_2(z_1) := \pm \frac{\sqrt{2}}{2}\sqrt{2d_2 - z_1^2 - \sqrt{(2d_2 - z_1^2)^2 - 4d_1^2}}, \]
$d_2 > d_1$ are positive real numbers.

By definition, it is easy to see that
\[ (c_1(z_1) + c_2(z_1))^2 = 2d_2 - z_1^2 \pm 4d_1^2, \]
\[ (c_1(z_1) - c_2(z_1))^2 = 2d_2 - z_1^2 \mp 4d_1^2. \]
Thus here we need $d_2 \geq 2 d_1^2$. Next we  prove $F$ is a special case of (\ref{symK=-1F}).

By a direct computation, we have
\[ c_1(z_1) + c_2(z_1) + z_2 = \sqrt{ 2d_2 \pm 4d_1^2 - |x|^2 + \frac{\langle x,y\rangle^2}{|y|^2}} +  \frac{\langle x,y\rangle}{|y|},  \]
\[ c_1(z_1) - c_2(z_1) - z_2 =\sqrt{ 2d_2 \mp 4d_1^2 - |x|^2 + \frac{\langle x,y\rangle^2}{|y|^2}} -  \frac{\langle x,y\rangle}{|y|}.  \]
\begin{equation*}
\begin{split}
F&   = \frac{|y| c_1(z_1)}{ c_1(z_1)^2 - (z_2 + c_2(z_1)^2)^2}   \\
&= \frac{|y|}{2} \frac{ c_1(z_1) - c_2(z_1) - z_2 +  c_1(z_1) + c_2(z_1) + z_2}{ ( c_1(z_1) - c_2(z_1) - z_2 )( c_1(z_1) + c_2(z_1) + z_2 )  }\\
   & = \frac{|y|}{2} \Big(  \frac{1}{c_1(z_1) + c_2(z_1) + z_2} + \frac{1}{ c_1(z_1) - c_2(z_1) - z_2} \Big)\\
   & = \frac{1}{2} \Big\{ \frac{\sqrt{(2d_2 \pm 4d_1^2 - |x|^2)|y|^2 + \langle x,y\rangle^2} - \langle x,y\rangle}{2d_2 \pm 4d_1^2 - |x|^2}  \\
   &+ \frac{\sqrt{(2d_2 \mp 4d_1^2 - |x|^2)|y|^2 + \langle x,y\rangle^2} +\langle x,y\rangle}{2d_2 \mp 4d_1^2 - |x|^2}  \Big\}.
\end{split}
\end{equation*}
Thus, $F$ is a  special case of (\ref{symK=-1F}) by a constant scaling.
\end{ex}
When ${\bf K}=1$, by Theorem \ref{thmK=1} we obtain
\begin{thm}
Let $F = |y|\zeta(|x|,\frac{\langle x,y\rangle}{|y|})$ be a spherically symmetric Finsler metric on on a convex domain $\Omega\subseteq R^n$.
Suppose $F(x,y)$ can be extended to a complex function $F(x,y+ xz)$, $z \in \mathbb{C}$.
Then $F$ is locally projectively flat with constant flag curvature
${\bf K} =1$ if and only if
\be \label{symK=1F}
F = \mathcal{I}m\Big[  \frac{(c+i)^2\langle x,y\rangle + \sqrt{ (c+i)^2(1-(c+i)^2 |x|^2)|y|^2+ (c+i)^4 \langle x,y\rangle^2}}{1-(c+i)^2 |x|^2} \Big],
\ee
where $c$ is a constant. Obviously, when $c=0$
\[ F_{c=0} = \frac{ \sqrt{ (1 + |x|^2)|y|^2 - \langle x,y\rangle^2}}{1+ |x|^2}. \]
\end{thm}
{\bf Proof}: The sufficiency can be verified directly. We only need to prove the necessity. By Lemma \ref{lemsym},
\[ F(0,y) =  |y|,\ \ \  P(0,y) = c |y|, \]
where $c$ is a constant.
Then by Theorem \ref{thmK=1}
\be
 \Psi =  c|y + x \Psi| + i |y+x \Psi|.
\ee
Solving above equation we get
\[ \Psi =  \frac{(c+i)^2\langle x,y\rangle + \sqrt{ (c+i)^2(1-(c+i)^2 |x|^2)|y|^2+ (c+i)^4 \langle x,y\rangle^2}}{1-(c+i)^2 |x|^2}. \]
Thus by (\ref{H}) we obtain (\ref{symK=1F}).
\qed
\begin{equation*}
\end{equation*}

{}
\noindent
Benling Li\\
Department of Mathematics\\
Ningbo University\\
Ningbo, Zhejiang Province 315211\\
P.R. China\\
libenling@nbu.edu.cn


\begin{thebibliography}{Mat2}

\bibitem{Be1} L. Berwald,  {\it Parallel$\acute{u}bertragung$ in allgemeinen R$\ddot{a}umen$}.  Atti
Congr. Intern. Mat Bologna \textbf{4} (1928), 263-270.
\bibitem{Be2} L. Berwald,  {\it $\ddot{U}$ber die n-dimensionalen  Geometrien konstanter Kr$\ddot{u}mmung$, in denen die Geraden die $k\ddot{u}rzesten$ sind,} Math. Z.\textbf{ 30}(1929), 449-469.

\bibitem{3} R. Bryant, {\it Finsler structures on the 2-sphere satisfying $K = 1$}, Finsler Geometry,
Contemporary Mathematics \textbf{196}, Amer. Math. Soc., Providence,
RI, 1996, 27-42. MR \textbf{97e}:53128.
\bibitem{4} R. Bryant, {\it Projectively flat Finsler 2-spheres of constant curvature}, Selecta Math.,
New Series, \textbf{3}(1997), 161-204. MR \textbf{98i}:53101.

\bibitem{F1} P. Funk, {\it $\ddot{U}$ber Geometrien, bei denen die Geraden die K$\ddot{u}$rzesten
sind}, Math. Annalen \textbf{101}(1929), 226-237.
\bibitem{F2} P. Funk, {\it $\ddot{U}$ber zweidimensionale Finslersche R$\ddot{a}$ume, insbesondere
$\ddot{u}$ber solche mit geradlinigen Extremalen und positiver konstanter
Kr$\ddot{u}$mmung}, Math. Zeitschr. \textbf{40}(1936), 86-93.

\bibitem{5} G. Hamel, {\it $\ddot{U}$ber die Geometrien in denen die Geraden die $K\ddot{u}rtzesten$ sind}, Math.
Ann. \textbf{57}(1903), 231-264.

\bibitem{Mo} L. Huang and X. Mo, {\it Projectively flat Finsler metrics with orthogonal invariance}, Annales Polonici Mathematici {\bf 107} (2013), 259-270.
\bibitem{LiShen} B. Li and Z. Shen,{\it On a class of projectively flat Finsler metrics with constant flag curvature}, Inter. Jour. Math, Vol.18, No. 7 (2007) 749-760.

\bibitem{Sh1} Z. Shen, {\it Projectively  Flat Finsler metrics of Constant Flag curvature}, Trans of Amer. Math. Soc. {\bf
355}(4)(2003), 1713-11728.

\bibitem{XuLi} B. Xu and B. Li, {\it On a class of projectively flat Finsler metrics with flag curvature $K=1$}, Diff. Geo. Appl.  {\bf 31}(2013), 524-532.

\bibitem{Zhou} L. Zhou, {\it Projective spherically symmetric Finsler metrics
with constant flag curvature in $R^n$}, Geom Dedicata {\bf 158}(2012), 353-364.


\end{thebibliography}
\end{document}